\documentclass[journal, 9pt]{IEEEtran}
\IEEEoverridecommandlockouts
\usepackage{cite}
\usepackage{amsmath,amssymb,amsfonts}
\usepackage{algorithmic}
\usepackage{graphicx, tabularx}
\usepackage{textcomp}
\usepackage{xcolor}
\usepackage{enumerate}
\usepackage{caption, subcaption}
\usepackage[compact]{titlesec}
\titlespacing{\section}{0pt}{0pt}{0pt}
    
\begin{document}
\pagestyle{plain}
\setlength{\abovedisplayskip}{5pt}
\setlength{\belowdisplayskip}{5pt}

\title{Impact of Electric Vehicle Routing with Stochastic Demand on Grid Operation}


\author{\IEEEauthorblockN{Oluwaseun Oladimeji, Alvaro Gonzalez-Castellanos, David Pozo,}
\IEEEauthorblockA{\textit{Skolkovo Institute of Science and Technology}\\ }
\and
\IEEEauthorblockN{Yury Dvorkin, Samrat Acharya,}
\IEEEauthorblockA{\textit{New York University}\\}
\thanks{
\hrule \vspace{.1cm}
This work was partially funded by the Skolkovo Institute of Science and Technology as a part of the Skoltech NGP Program (Skoltech-MIT joint project).}
\vspace{-0.6cm}
}

\maketitle
\IEEEpubidadjcol

\vspace{-0.5cm}

\begin{abstract}
Given the rise of electric vehicle (EV) adoption, supported by government policies and dropping technology prices, new challenges arise in the modeling and operation of electric transportation.
In this paper, we present a model for solving the EV routing problem while accounting for real-life stochastic demand behavior.
We present a mathematical formulation that minimizes travel time and energy costs of a EV fleet. The EV is represented by a battery energy consumption model. To adapt our formulation to real-life scenarios, customer pick-ups and drop-offs were modeled as stochastic parameters. A chance-constrained optimization model is proposed for addressing pick-ups and drop-offs uncertainties. Computational validation of the model is provided based on representative transportation scenarios. Results obtained showed a quick convergence of our model with verifiable solutions. Finally, the impact of electric vehicles charging is validated in Downtown Manhattan, New York by assessing the effect on the distribution grid.
\end{abstract}


\section*{Nomenclature}
\addcontentsline{toc}{section}{Nomenclature}
\begin{IEEEdescription}[\IEEEusemathlabelsep\IEEEsetlabelwidth{$i,j,s \in J_0$}]
\item[\textit{Indexes and Sets}]
\item[$i,j,s \in J$] Customer locations.
\item[$i,j,s \in J_0$] Transportation nodes, i.e., customers locations and depot; $|J_0|=|J|+1$.
\item[$k \in K$] Charging stations. $K = J_0 \cap N$.
\item[$n, m \in N$] Nodes in the electric distribution network.
\item[$(n,m) \in \mathcal{L}$] Power distribution  lines.
\item[$v \in V$] Electric vehicles.
\item

\item[\textit{Parameters}]
\item[$\epsilon_j$] Risk tolerance for demand satisfaction at node \(j\),  \hfill[-].
\item[$C_v$] Vehicle power consumption rate, \hfill[kWh/min].
\item[$C^\text{E}$] Cost of electricity,  \hfill[\$/kW].
\item[$C^\text{T}$] Unitary cost of time travelled,  \hfill[\$/min].
\item[$D_{i}$] Amount of dropped off customers at $i$, \hfill[passengers].
\item[$\underline{E}_v/\overline{E}_v$] Minimum/maximum battery charge level,  \hfill[kWh].
\item[$ET_{i}$] Earliest time of arrival at $j$,  \hfill[min].
\item[$\overline{L}_v$] Vehicle capacity, \hfill[passengers].
\item[$LT_{i}$] Latest time of arrival at $j$,  \hfill[min].
\item[$M^{(\cdot)}$] Large constant,  \hfill[--].
\item[$P_{i}$] Amount of picked up customers at $i$,  \hfill[passengers].
\item[$\underline{P}_n/\overline{P}_n$] Minimum/maximum nodal active demand,  \hfill[kW].
\item[$P^\text{de}_n/Q^\text{de}_n$] Local active/reactive demand at node $n$,  \hfill[kW/kVar].
\item[$\underline{Q}_n/\overline{Q}_n$] Minimum/maximum nodal reactive demand,  \hfill[kVar].
\item[$R_v$] Vehicle charging rate,  \hfill[kWh/min].
\item[$\underline{S}_n/\overline{S}_n$] Minimum/maximum nodal apparent demand,  \hfill[kVA].
\item[$T_{ij}$] Travel time between nodes,  \hfill[min].
\item[$\overline{T}_k/\underline{T}_k$] Minimum/maximum charging time,  \hfill[min].
\item[$\overline{T}^\text{rech}_v$] Maximum allowed time before recharge,  \hfill[min].
\item[$\underline{V}/\overline{V}$] Minimum/maximum nodal voltage,  \hfill[V].
\item

\item[\textit{Variables}]
\item[$\pi_j$] Variable used to prohibit sub-tours,  \hfill[--].
\item[$\tau_{kv}$] Charging time spent by vehicle $v$,  \hfill[min].
\item[$e_{jv}$] Energy level of vehicle $v$ at node $j$,  \hfill[kWh].
\item[$l^\text{0}_v$] Vehicle load after leaving the depot,  \hfill[passengers].
\item[$l_{jv}$] Vehicle load after visiting location $j$,  \hfill[passengers].
\item[$p_{nm}/q_{nm}$] Active/reactive power flow along line $\left(n,m\right)$,  \hfill[kW/kvar].
\item[$p^\text{d}_n$]Sum of electric power consumption and EV charging power at node $n$,  \hfill[kW/kvar].
\item[$p^\text{g}_n/q^\text{g}_n$] Active/reactive power generation at $n$,  \hfill[kW/kVar].
\item[$t_{jv}$] Arrival time of $v$ to node $j$,  \hfill[min].
\item[$u_n$] Square of the voltage magnitude at $n$,  \hfill[$V^2$].
\item[$w_{kv}$] Auxiliary variable combining charge duration and power usage at station, \hfill [kWh].
\item[$x_{ijv}$] Binary variable indicating whether vehicle $v$ travels directly from $i$ to $j$,  \hfill[--].
\item[$y_{kv}$] Binary variable indicating vehicle $v$ at charging station $k$,  \hfill[--].
\end{IEEEdescription}


\section{Introduction} \label{introduction}
Until recently, large scale electrification of the transportation sector has not gained so much attention, even though transportation represents more than 25\% of the world's energy consumption \cite{ieaconstransportation}. The transportation sector contributes to a large share of greenhouse emissions, e.g., in California, it accounts for more than 40\% of emissions \cite{californiaGHG}. Energy decarbonization road-maps highlight electrification of the transport sector as a critical step to reach emission goals. Likewise, emissions reductions are especially necessary for polluted mega-cities, where an increasing number of projects have replaced combustion-based bus fleets by electric ones.
%
This modernization benefits from improvements in rechargeable battery technology leading to quantitative improvement in EV technology. EV manufacturing companies continue to make measurable production growth with a significant presence of EVs used for commercial transportation services \cite{globaloutlook}.

In terms of impact of vehicles charging on the distribution grid, let us evaluate how much power is drawn from the grid during routing. Using a \mbox{Tesla model S} with a battery capacity of 75kWh with consumption of only about 33kWh for 100 miles driven translates to a driving range of about 220 miles. This might mean only a single charge in 2 days for private car owners which sit idle for about 95\% of the time on a daily basis \cite{cars95}. However, when we consider that ride sharing where cars would constantly be on the road, then the effect on the electric network are: more frequent charging sessions and increasing use of fast-charging stations.

If we consider fast chargers with power charging rate of 50kW, twenty vehicles would be withdrawing 1MW from the power network if they are charged simultaneously.
New York uses about 11000MWh of electricity per day which translates into an average hourly consumption of approximately 450MW. Thus, only twenty vehicles would represent about 0.25\% of the city's demand.
An average ride-sharing car makes between 5 to 12 trips per day. Thus, an EV would require charging at least once or twice per day. If we consider a fleet owner with 500 EVs, power requirement increases up to 25MW 
a considerable increase on demand is requested to the electric network.

Transportation electrification then calls for modification of previous vehicle routing algorithms, which have considered traditional internal combustion engine (ICE) vehicles. ICE vehicles estimate energy consumption by distance traveled and gross weight of the vehicle. For EVs, these assumptions are still valid, but additional considerations such as battery properties and charging infrastructure must also be considered.

\subsection{Literature Survey} \label{literature}
The vehicle routing problem (VRP) with simultaneous pickup and delivery for material goods was introduced in \cite{min1989multiple}. In \cite{dethloff2001vehicle}, customer satisfaction was added as a model constraint. Time-window constraints are included in the VRP to prevent situations in which customers wait for long periods. A review on exact algorithms for time-and- capacity-constrained models was presented in \cite{baldacci2012recent}.
Including uncertainty, \cite{probcapacitymahdi} modelled the stochastic nature of drop-offs and pickups using chance constraints. In another study, non-linear interval-based programming was used to investigate VRP with uncertainty in drop-offs \cite{cao2018research}.

As energy management became increasingly important, various studies incorporated partial battery charging stations into their model while still using a conventional combustion engine model for energy management. The green-mixed fleet VRP was proposed in \cite{macrinaguisy}. Models with minimization of traveled distance and energy usage have also been presented in \cite{hiermann2016electric, schneider2015adaptive, zhang2018electric, ferro2018optimization}.

Single user vehicles, due to their passive usage, are primarily considered as static energy storage equipment. However, in this paper we seek to address the question on EVs applicability based on two premises: 1) EVs production continues to increase, and 2) ride/car sharing popularity, which increases vehicle mileage. As EVs increase in numbers, there will be an increasing request for power reliability from the existing power infrastructure. Additionally, if EVs are widely accepted as an alternative to commercial vehicles, we pose a further question of how much impact this would have on the power grid.

The current solution being proposed to curtail the impact of EVs on the power grid is valley filling. Inherent challenges with this is the possibility of creating a new peak and importantly, controlling human behaviour. Valley filling \cite{jian2017high} essentially seeks to smooth out the individual or aggregate load curve by shaving peak loads and increasing off-peak loads e.g. charging EVs at night.

In \cite{chenoptimal,chenvalley}, offline and online valley filling algorithms were introduced. The authors recognised the challenge posed by large fleet of EVs on the distribution network if the EVs are stationed in residential areas. The EV charging problem was modelled using a modified version of the optimal power flow problem with a valley-filling objective. However, electricity price was considered constant in the papers and EVs were taken as stationary.
A mode of decentralized coordination of charging large populations of EVs using non-cooperative game theory was developed in \cite{callawaydecen}. The charging schedule was also modeled as a valley-filling problem. As in the previous literature discussed, electricity price was taken as constant and uncertainties were not considered.
A control algorithm focusing on peak shaving was proposed to match the power curve and scheduled EV usage curve in \cite{wangv2g}. A grid that incorporates smart meters, intelligent switches, two-way charge devices and a central control center was considered. The charging decision is then made by matching the load forecast curve and the available charge/discharge curve of EVs.

A major drawback in these systems is the absence of mobility-awareness when scheduling EV charging. An additional challenge was lack of real data in justifying the proposed control algorithms/frameworks developed. This is also a reason why authors modeled electricity prices as time-invariant due to lack of evidence regarding how the EVs will affect the dynamics of the power system.

Furthermore, several tools have been developed to address some of the questions that we pose in our work. CASPOC, COMPOSE and HOMER are three simulators that are closest to our proposed methodology \cite{mahmud2016review}. However, these tools have limitations which do not capture the solution that we propose. HOMER, for instance, does not include mobility and time dependent batteries based on specific EVs could not be modeled on it. CASPOC and COMPOSE, though can model EVs, do not support integration of mobility with the grid.

This paper introduces a new variant of co-optimization of multi-energy systems which have included coordinated operation of power, thermal and gas networks. An operational optimization model for electricity, thermal and gas  network was presented in \cite{wang2019operation}. Coordinated operation of the heating and electricity system for the United Kingdom is discussed in \cite{zhang2018economic}. A survey on power and gas system integration is presented \cite{farrokhifar2020energy}. These coordinated systems offer environmental, system reliability, and economic benefits; similar benefits can be achieved by the coordination of the transportation and electricity network.

\subsection{Paper Contribution and Organization}
The main contribution of this work is to propose a co-optimization model that combines the routing of an EV transportation fleet and the distribution grid power flow. Additionally, we extend the proposed model by incorporating pick-up and drop-off uncertainty in the VRP model. The resulting model is formulated as a chance-constrained optimization problem and reformulated as an equivalent deterministic mixed-integer linear programming (MILP). Then, we present an illustrative case study based on real data of the island of Manhattan, New York.

The rest of the paper is organized as follows. Section \ref{Maths model} develops the mathematical formulation of the VRP problem with vehicle recharging, demand uncertainty and power flow in the electric distribution network. In Section \ref{Case study}, we describe the test case data, and how it is processed and used for computational tests. In Section \ref{Results} the proposed model is tested on the case study. Section \ref{Conclusion} concludes the study.


\section{Mathematical Model} \label{Maths model}
The overall aim of the optimization problem is to reduce the EV operation cost by minimizing travel time and energy usage of the vehicles and to minimize the cost of the distribution network operation by efficiently scheduling EV charging periods. This section is devoted to formulating electric transport routing, representative energy management model and the distribution grid. Noting that some parameters are not exactly known during planning, equations used to describe stochasticity will also be highlighted. 
Our goal is to provide a routing and EV charging plan for EV-fleet owners and electric utility companies. This plan will enable a reduction in distance traveled, contribute to reduced energy usage and consequently lead to efficient and operation of the distribution grid.

\subsection{Electric Vehicle Routing Model} \label{Routing model}
The vehicle routing problem consists of assigning routes to a fleet of transportation vehicles based on drop-off and pick-up needs of customers distributed on a set of nodes on the transportation network. The objective \eqref{eq:distobjective} of the VRP is the minimization of total time traveled by the transportation fleet.
$x_{ijv}$ represents whether the vehicle $v$ travels directly from stop $i$ to stop $j$, and $T_{ij}$ is the distance between these stops.
In combination with \eqref{eq:distobjective}, the VRP can be modelled by the following problem.
%
\begin{subequations} \label{eq:VRP}
\begin{IEEEeqnarray}{lll}
    \min & \sum\limits_{i \in J_0}\sum\limits_{j \in J_0}\sum\limits_{v \in V} x_{ijv}T_{ij}  \label{eq:distobjective} \\
    \text{s.t.: } \:
    &
    \sum\limits_{i \in J_0}\sum\limits_{v \in V} x_{ijv} = 1,  & \forall{j} \label{eq:1 vehicle serve one customer}\\ 
    &\sum\limits_{i \in J_0} x_{isv} = \sum\limits_{j \in J_0} x_{sjv}, & \forall{s \in J, v} \label{eq:served by 1 vehicle} \\
    &l^0_v = \sum\limits_{i \in J_0}\sum\limits_{j \in J} D_jx_{ijv}, & \forall{v} \label{eq:initial vehicle load from the depot} \\
    &l_{jv} \geq l_{iv} - D_j + P_j- M^L_{ijv}, & \forall{i\ne j,v} \IEEEeqnarraynumspace\label{eq:vehicle capacity conservation on tour}\\
    &l_{jv} \leq \overline{L}_v, & \forall{j,v} \label{eq:cap of vehicle on tour}\\
    &\pi_j \geq \pi_i + 1 - |J_0|(1- \sum\limits_{v \in V} x_{ijv}), & \forall{i\ne j} \label{eq:sub tour breaking} \\
    &ET_j \leq t_{jv} \leq LT_j, & \forall{j, v} \label{eq:bound on arrival time} \\
    &t_{jv} \geq t_{iv} + \tau_{kv} + x_{ijv} T_{ij} - M^T_{ijv}, \qquad & \forall{i {\in}  J, j {\in} J_0, k, v} \IEEEeqnarraynumspace \label{eq:arrival time eqn} \\
    &l_{jv}, t_{jv}, \pi_j \geq 0,&  \forall{j,v} \label{eq:sub-tour breaking bound}\\
    &x_{ijv} \in \{0,1\}, & \forall{i,j,v}. \label{eq:x_nature}
\end{IEEEeqnarray}
\end{subequations}
Constraints \eqref{eq:1 vehicle serve one customer} and \eqref{eq:served by 1 vehicle} define that each customer should be moved from a pick-up point to a drop-off point by exactly one vehicle.
Initial vehicle load from the depot is set by  \eqref{eq:initial vehicle load from the depot}. \eqref{eq:vehicle capacity conservation on tour} represents the vehicle loads en route; with $M^L_{ijv} = M^L(1-x_{ijv})$. The constant $M^L$ is a large number that sets a common lower bound to constraints in which $x_{ijv}$ has a unitary value. The value of $M^L$ can be calculated beforehand based on demand levels of the transportation network.
Vehicle capacity limits are set by \eqref{eq:cap of vehicle on tour}. 

As the algorithm searches for solutions to satisfy given constraints, many disjointed routes which have no connections with the depot will be created. However, a fleet manager requires all vehicles to be returned to the starting station at the end of every tour. Sub-tour eliminating constraint, \eqref{eq:sub tour breaking}, removes solutions without connections to the depot.

Based on customer-established pick-up times, expression \eqref{eq:bound on arrival time} sets the bounds for the arrival time at different stops. Arrival time is calculated with \eqref{eq:arrival time eqn}, where the value $M^T_{ijv}$ is calculated analogously as $M_{ijv}^L$, but with a large time constant $M^T$.
The lower bounds on vehicle load, traveled time and sub-tour breaking variable is set by \eqref{eq:sub-tour breaking bound}, whereas \eqref{eq:x_nature} defines $x_{ijv}$ as a binary variable.

\subsection{Energy Usage Model} \label{Energy model}
For the routing of the electric fleet, it is necessary to model energy usage of the vehicle and time spent at charging stations. The objective of the VRP with energy management can be set either as the minimization of the total time spent at charging stations or as the total cost of energy used. If a constant electricity price is considered, both objectives are equivalent. Thus, we express the objective of the VRP with energy management as:
\begin{subequations}
\begin{IEEEeqnarray}{lll} 
 \min & \sum\limits_{k \in K}\sum\limits_{v \in V}\tau_{kv},   \label{eq:Energyobjective} \\
\text{s.t.: } \:
    & e_{jv} \leq e_{iv} {+} \tau_{kv}R_v {-} C_v x_{ijv}T_{ij} {+} M^E_{ijv}{,} \quad  & \forall{i {\ne} j, k, v} \label{eq:Energy conservation during route} \\
    &\sum\limits_{v \in V}x_{ijv} T_{ij} \leq \overline{T}^\text{rech}, & \forall{i, j} \label{eq:max travel time of 1 vehicle} \\
    &\underline{T}_k y_{kv}\leq \tau_{kv} \leq y_{kv}\overline{T}_k, & \forall{k, v} \label{eq:bounds on charging time}\\
    &e_{kv} + \tau_{kv}R_v \leq \overline{E}_v, & \forall{k, v} \label{eq:bound on max battery energy during charging}\\
    &\underline{E}_v \le e_{jv} \leq \overline{E}_v, & \forall{j, v}. \label{eq:bound on battery energy levels}
\end{IEEEeqnarray}
\end{subequations}

Battery energy levels at each node are calculated by \eqref{eq:Energy conservation during route}. $M_{ijv}^E = M^E(1-x_{ijv})$ is a large value that is used to activate this constraint \eqref{eq:Energy conservation during route} when the edge connecting the stops $i$ and $j$ is traveled.
Expression \eqref{eq:max travel time of 1 vehicle} sets an upper bound on the time a vehicle can move before returning to depot. 
\eqref{eq:bounds on charging time} and \eqref{eq:bound on max battery energy during charging} respectively define the lower and upper bounds on the time of charging for a vehicle visiting a charging station and the energy after a charging stop. Battery capacity limits are defined by \eqref{eq:bound on battery energy levels}.
Battery degradation is not considered in our formulation due to modern battery technology which are designed to ensure a good lifetime. Studies have shown that EV battery degradation rate is non-linearly correlated with its usage, and for EV fleets, increased utilization does not necessarily impact battery life \cite{geotabevba3health}.

\subsection{Distribution Grid Model} \label{grid model}
The conventional methods for solving power flow in distribution grid follow the fundamental electrical current flow laws. We formulate the electrical distribution grid using LinDistFlow  equations \cite{baran1989optimal} which is a linear form of the power flow equations in electrical distribution grids, resulting into a  computationally and numerically efficient formulation for distribution network operation modeling.
\begin{subequations} \label{eq:grid model}
\begin{IEEEeqnarray}{lCll}
    p_{nm} &=&  p_m^\text{g} - P_m^\text{de} - \sum_{m':(m,m') \in \mathcal{L}} p_{mm'}, \quad  & \forall{(n,m) \in \mathcal {L}}  \label{eq:activepowerdiff} \\
    q_{nm} &=&  q_m^\text{g} - Q_m^\text{de} - \sum_{m':(m,m')  \in \mathcal{L}} q_{mm'}, \quad  & \forall{(n,m) \in \mathcal {L}}  \label{eq:reactivepowerdiff} \\
    u_m &=& u_{n} - 2(r_{nm}p_{nm}+x_{nm}q_{nm}), \qquad & \forall{(n,m) \in \mathcal {L}} \IEEEeqnarraynumspace \label{eq:voltage} \\
    p_{nm} &\leq& \overline{P}_{nm}, \quad & \forall{(n,m) \in \mathcal {L}} \label{eq:activepowergen} \\
    q_{mn} &\leq& \overline{Q}_{nm}, \quad & \forall{(n,m) \in \mathcal {L}}  \label{eq:reactivepowergen}\\
    \underline{P}_n &\leq& p_n^\text{g} \leq \overline{P}_n, \quad & \forall{n}  \label{eq:activepowercons} \\
    \underline{Q}_n &\leq& q_n^\text{g} \leq \overline{Q}_n, \quad & \forall{n}  \label{eq:reactivepowercons} \\
    \underline{V}_n^2 &\leq& u_n \leq \overline{V}_n^2, \quad & \forall{n}, \label{eq:voltagecons}
\end{IEEEeqnarray}
\end{subequations}
\noindent Where $P^\text{de}_n + iQ^\text{de}_n$ represents the complex power consumption of electric demand only at node $n$. $p_{nm} + iq_{nm}$ is complex power flowing along the power line $\left(nm\right)$. $p_n^\text{g} + iq_n^\text{g}$ is the power generation at node n. $u_n$ is the square of the nodal voltage. Equations \eqref{eq:activepowerdiff} and \eqref{eq:reactivepowerdiff} model power flow through the lines. The nodal voltage is calculated by \eqref{eq:voltage}; while line flow, generation, and voltage limits are set by \eqref{eq:activepowergen}--\eqref{eq:voltagecons}.

\subsection{Transportation Network and Electrical Grid Coupling}
The transportation network and the electric grid are coupled by EV charging stations. The power required by EVs during charging should be added to the grid's existing demand. 
Note that the location and power demand from EVs are dependent on the routes resultant of the VRP.
Electricity consumption in the distribution grid can be summed up to that in the transportation grid at the nodes where charging stations are present by the following expression:
\begin{IEEEeqnarray}{r"l}
p^\text{d}_n =  p^{\text{d},e} + R_v,  & (\forall  n \in J, \forall k \in K) \wedge y_{kv} = 1. \label{eq:coupling_eq}
\end{IEEEeqnarray}
Thus, active power consumption $p^\text{d}_n$ at node $n$, is equivalent to regular electrical load $p^{\text{d},e}$, plus electric load related with transportation grid, $R_v$,  at charging station $k$. 

\subsection{Modelling Uncertainty} \label{uncertainty}
Given the uncertain nature of demand, we need to express the uncertainty in the drop-off and pick-up parameters of the transportation demand balance \eqref{eq:vehicle capacity conservation on tour}. We opt for the use of chance-constrained programming to describe demand uncertainty because this approach allows us to satisfy customer demands with a certain degree of confidence, i.e., we model the constraints in such a way that drop-offs and pick-ups are complied with a probability above the established level ($1-\epsilon_j$):
\begin{IEEEeqnarray}{llr}
    \mathbb{P}(l_{jv} \geq l_{iv} - D_j + P_j - M^L_{ijv}) \geq 1 - \epsilon_j, & \quad \forall{i\ne j},v \label{eq:vehicle capacity conservation when leaving depot 2}
\end{IEEEeqnarray}
where $l_{iv}$ and $l_{jv}$ are variables representing vehicle nodal loads, and $D_j$ and $P_j$ are bounded parameters. The probability of the nodal constraint satisfaction is given by $(1 - \epsilon_j) \in [0,1]$. The use of nodal constraint satisfaction allows us to assign differential priority to transportation nodes based on their importance by selecting different values of the risk parameter $\epsilon_j$, i.e., more important transportation nodes have lower risk values \(\epsilon_j\). The following reformulation of vehicle capacity conservation is made \cite{lubin2019chance}:
\begin{IEEEeqnarray}{lll}
    l_{jv} {\geq} l_{iv} {-} \mathbb{E}[D_j-P_j] \pm \Phi^{-1}_{1{-}\epsilon_j} \sigma_j {-} M^L_{ijv}, &\quad \forall{i\ne j},v \IEEEeqnarraynumspace\label{eq:vehicle capacity conservation on tour 3}
\end{IEEEeqnarray}

The value $\mathbb{E}[D_j-P_j]$ represents the net demand expected at node $j$. $\Phi^{-1}_{1{-}\epsilon_j}$ is the inverse cumulative distribution function evaluated at $1{-}\epsilon_j$, and $\sigma_j$ the standard deviation of the nodal net demand.

\subsection{Overall Objective} \label{overall objective}
We consider writing a single objective by making some adjustments to our prior objectives defined above.
\begin{itemize}
    \item A penalization factor $C^\text{T}$ is introduced to assign a cost related to traveling time in \eqref{eq:distobjective}. This factor is similar to a transport service cost based on distance or time traveled.
    \item By considering power consumption and time spent charging at different nodes, multiplied by electricity cost, it is possible to obtain the total cost of energy. The electricity cost is constant here because of the resolution of our planning horizon which is within an hour or less and the electricity cost is generally non-variant within this time period. Thus, the electricity cost will be updated according to prevalent prices during the time when the planning schedule is required.
\end{itemize} The new objective function becomes:
\begin{IEEEeqnarray}{C}
    \min \:  \sum\limits_{i,j,v}C^\text{T} x_{ijv}T_{ij} +  \sum\limits_{v,k} C^\text{E} p_k^c \tau_{kv} \label{eq:overall objective}
\end{IEEEeqnarray}

The second part of the final objective \eqref{eq:overall objective} introduces non-convexity to the model because both power consumption and charging time are variables. However, power consumption and charging time are both bounded by battery capacity and charge rate. We overcome the non-convexity by combining these two variables via the introduction of a McCormick envelope \cite{mccormick} which is represented by the set of constraints \eqref{eq:mccormick}: 
\begin{subequations}  \label{eq:mccormick}
\begin{IEEEeqnarray}{llr}
    w_{kv} \leq \overline{E}_v\tau_{kv} + p_k^c\underline{T}_k - \overline{E}_v\underline{T}_k,\qquad  & \forall{k, v} \label{eq:energy&time1} \\
    w_{kv} \leq p_k^c\overline{T}_k + \underline{E}_v\tau_{kv} - \underline{E}_v\overline{T}_k,\qquad  & \forall{k, v} \label{eq:energy&time2} \\
    \underline{E}_v \leq p_k^c \leq \overline{E}_v, \\
    \underline{T}_k \leq \tau_{kv} \leq \overline{T}_k,\qquad & \forall{k, v} \label{eq:tau}
\end{IEEEeqnarray}
\end{subequations}

A linear representation of the objective is finally given by:
\begin{IEEEeqnarray}{C}
    \min \sum\limits_{i,j,v} C^\text{T}x_{ijv}T_{ij} + \sum\limits_{v,k}C^\text{E}w_{kv} \label{eq:overall objective2}
\end{IEEEeqnarray}

\subsection{Model Summary} \label{model summary}
For further analysis, we classify our formulation into two broad groups: deterministic and chance-constrained stochastic models. 
As shown in Table~\ref{table:model summary}, the main difference in formulation appears in the set of constraints defining the EV routing problem. 
\begin{table}[ht]
\begin{center}
\caption{Model Summary}
\begin{tabular}{l | r | r} 
 \hline\hline
   & \multicolumn{1}{c|}{Deterministic} & \multicolumn{1}{c}{Stochastic CC} \\ [0.5ex]
 \hline
    Objective  & \eqref{eq:overall objective2} & \eqref{eq:overall objective2}  \\     
   
    EV routing problem  & \eqref{eq:1 vehicle serve one customer} -- \eqref{eq:x_nature} & \eqref{eq:1 vehicle serve one customer} -- \eqref{eq:initial vehicle load from the depot}, \\
    
    &  & \eqref{eq:cap of vehicle on tour} -- \eqref{eq:x_nature}, \eqref{eq:vehicle capacity conservation when leaving depot 2} \& \eqref{eq:vehicle capacity conservation on tour 3}	\\ 
    Energy usage  & \eqref{eq:Energy conservation during route} -- \eqref{eq:bound on battery energy levels} & \eqref{eq:Energy conservation during route} -- \eqref{eq:bound on battery energy levels}	\\ 
    
    Electric distribution OPF  & \eqref{eq:grid model}  & \eqref{eq:grid model}	\\
    
    Coupling constraint  & \eqref{eq:coupling_eq}  & \eqref{eq:coupling_eq} \\
    
    Total charged energy  & \eqref{eq:mccormick}  & \eqref{eq:mccormick} \\
    \hline
 \hline
\end{tabular}
\label{table:model summary}
\end{center}
\end{table}


\section{Case Study} \label{Case study}
\subsection{Transportation Network Data}
The case study selected for validating our proposed methodology is based on the Downtown Manhattan transportation network. Drop-off and pick-up data is taken from taxi trip records by the New York City Taxi and Limousine Commission for 2018 \cite{nyctaxidata}. The data set includes pick-up and drop-off dates, times, locations, trip distances and driver-reported passenger counts. Table~\ref{table:data snippet} shows a snippet from this data set where PULID and DOLID are respectively the pick-up location ID and drop-off Location IDs, which indicate exact latitudes and longitudes of these locations.

\begin{table}
\begin{center}
\caption{Representative Data Snippet}
\begin{tabular}{c c r r r r} 
 \hline
  Pickup & Drop-off & \multicolumn{1}{c}{Trip} & \multicolumn{1}{c}{PULID} &  \multicolumn{1}{c}{DOLID} & \multicolumn{1}{c}{Passenger} \\ 
 date-time & date-time & \multicolumn{1}{c}{distance} &  &  & \multicolumn{1}{c}{count} \\ [0.5ex]
 \hline\hline
    12/31/2018  & 1/1/2019  & 16.6 & 162 & 26 & 1\\     
    11:58:45 & 12:37:04 & & & & \\
    \hline
    12/31/2018  & 1/1/2019 	& 3 & 144 & 170 & 2 \\     
    11:58:44 & 12:16:54 & & & & \\
    \hline
    12/31/2018 	& 1/1/2019 	& 6.4 & 162 & 13 & 1 \\
    11:58:37 & 12:14:04 & & & & \\
    [1ex]
 \hline
\end{tabular}
\label{table:data snippet}
\end{center}
\vspace{-0.5cm}
\end{table}
In order to use the data in our proposed VRP model, it was necessary to process through clusterization, while preserving its accuracy. The first opportunity for aggregation comes from passenger count. Taxis routinely pick up just one passenger per trip. Here, we extend our focus beyond one passenger scenario to ride-sharing cases. With ride-sharing, passengers can be added to each vehicle in the fleet until it reaches its capacity. This aggregation makes it possible to add all demand for an extended period, for example all drop-offs and pickups between 6am and 11am are aggregated and the case is evaluated such that all the demands are within a 2-hour interval.

Using latitudes and longitudes, we divided the data into distinct locations such as commercial, industrial, and residential regions. With this, we were able to track movement patterns at different times of the day. For example, people generally move from residential areas to commercial centers in the mornings and vice versa in the evenings. The routes from the solution of the model ensures that vehicles can pick up passengers until their capacity is exhausted while still fulfilling the pick-up time and other constraints.

\begin{figure}[ht]
     \centering
     \begin{subfigure}[t]{0.43\columnwidth}
         \centering
         \includegraphics[width=\columnwidth]{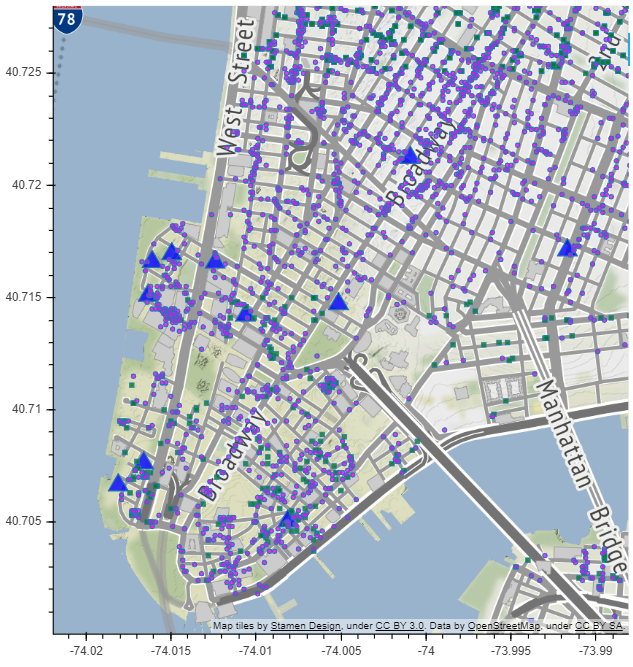}
         \caption{Passenger concentration and EV CS in Downtown Manhattan btw 6-10am, 13 Aug '18}
         \label{fig:Downtown manhattan}
     \end{subfigure}
     \hfill
     \begin{subfigure}[t]{0.43\columnwidth}
         \centering
         \includegraphics[width=\columnwidth]{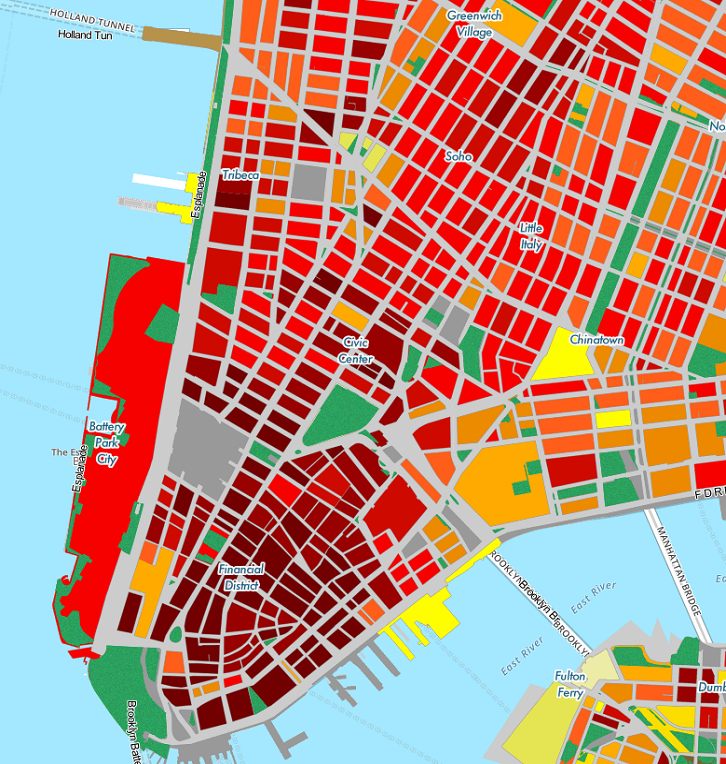}
         \caption{Downtown Manhattan energy consumption by area}
         \label{fig:Downtown manhattan power}
     \end{subfigure}
        \caption{Downtown Manhattan transport and energy demands}
        \label{fig:Downtown manhattan Transport and Energy Consumption}
\end{figure}

Fig.~\ref{fig:Downtown manhattan} shows a geographical snapshot of the demand data, which is at the Downtown Manhattan grid, where green squares represent originating nodes, red circles represent destination nodes and large blue triangles are EV charging stations.

\subsection{Electric Network Data}
For local electricity consumption, we obtained information from \cite{howard2012spatial} and New York City Mayor's Office of the Long-Term Planning and Sustainability \cite{nycSIRR}.
Electricity consumption at each block and lot level was obtained from \cite{howard2012spatial} and consumption patterns are shown by the different shades of red in Fig.~\ref{fig:Downtown manhattan power}. Darker red and lighter/yellow colors show regions of high and low power consumption respectively. Voltage ratings of each charging station location and reactance between stations were estimated similarly to \cite{acharya2020public}.




\section{Computational tests} \label{Results}
We compare the results of the deterministic and stochastic optimization here using the cases defined in Section \ref{model summary}.
For the deterministic model, known input include data outlined in Table~\ref{table:data snippet} and power requirements while variables to be determined include the route and energy consumption of each vehicle. For the stochastic model, demands at each node is unknown and modeled after a Gaussian distribution following analysis of historical data.

For clarity, we chose a 15-location case from the original data set with about 600 passengers over a 2-hour period. We consider a modest number of 6 vehicles deployed by the fleet owner to meet these demands, resulting in a model with over 1350 variables. As number of locations and vehicles increases, so does the dimension of $x_{ijv}$. Since there are 36 EV charging stations in Downtown Manhattan, (Fig.~\ref{fig:Downtown manhattan}), distance mapping is used to determine the nearest charging station when vehicles run out of charge or when power available would not be enough to reach the next passenger location.

\subsection{Routing Results}
\begin{figure}[ht]
     \centering
     \begin{subfigure}[b]{0.48\columnwidth}
         \centering
         \includegraphics[width=\columnwidth]{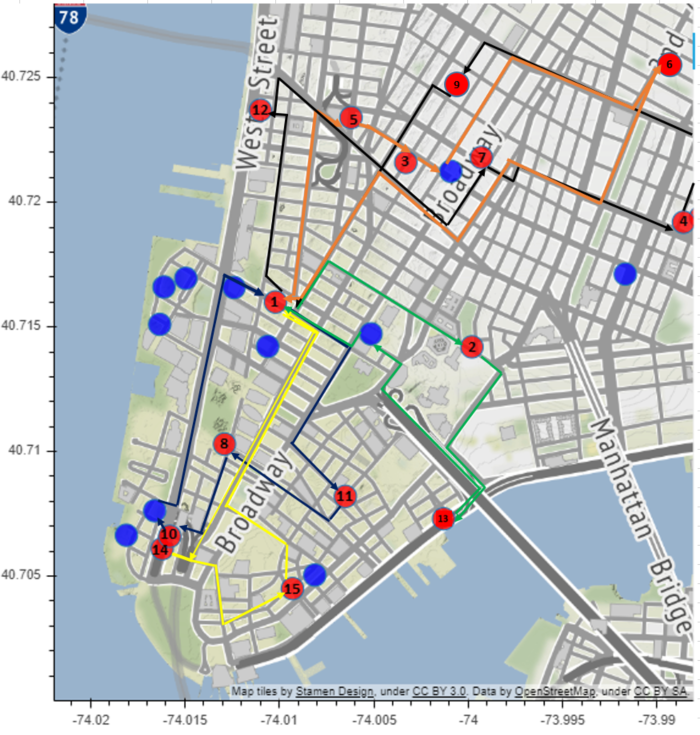}
         \caption{Routing plan for deterministic model}
         \label{fig:Routing plan_deterministic}
     \end{subfigure}
     \hfill
     \begin{subfigure}[b]{0.48\columnwidth}
         \centering
         \includegraphics[width=\columnwidth]{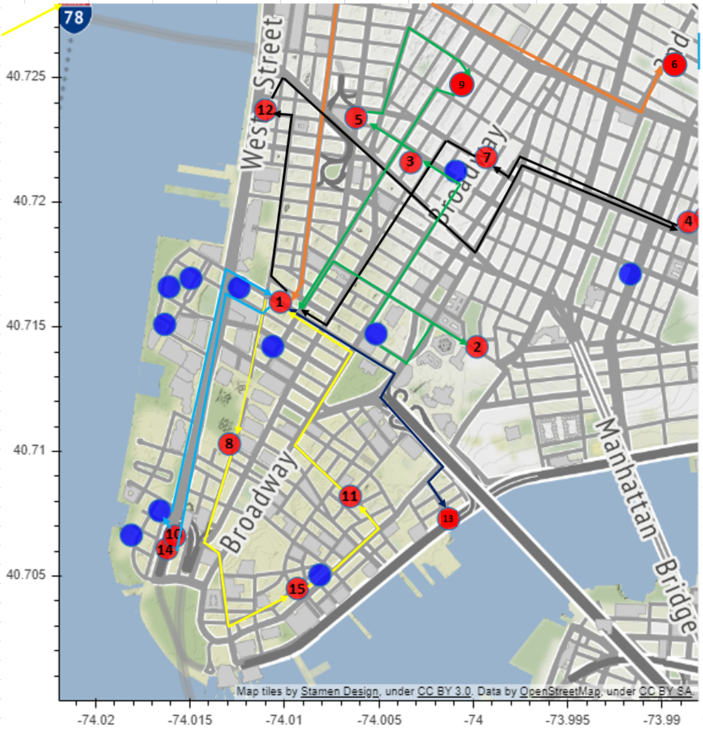}
         \caption{Routing plan for chance-constrained model}
         \label{fig:Routing plan with Chance Constraint}
     \end{subfigure}
        \caption{Results of vehicle routing problem}
        \label{fig: Vehicle Routing}
\end{figure}

Using historical data, Fig.~\ref{fig:Routing plan_deterministic} shows the routing results under assumption of perfect demand knowledge, i.e., deterministic modeling. Only 5 out of 6 available vehicles were deployed from the depot. However, when uncertainty is introduced in this scenario, result of routes varies from the deterministic model. It is shown in Fig.~\ref{fig:Routing plan with Chance Constraint} that the routing is changed and also an additional vehicle is deployed. This is the case when a high degree of confidence is used, i.e., $1-\epsilon_j \geq 99.5\%$.

However, for $90\% < 1- \epsilon_j < 95\%$, 5 vehicles were deployed, but with a different set of routes. The confidence levels translate to how many passengers will be picked or dropped in a particular location. Thus, as the tolerated risk varies, the optimal amount of vehicles and routes also change.

Table~\ref{table:other results} presents the optimization results of different test cases. Labels starting with 'D' and 'U' represent deterministic and stochastic problems respectively. Here, we contrast results of deterministic versus stochastic chance-constrained model. For few locations, $|J| \leq 25$, number of routes and objective value have similar value for both types of models. However, as number of locations $|J|$ increases, number of routes differ greatly between both models. The same trend is seen with the objective value. The number of routes is equivalent to number of vehicles deployed from the fleet while the objective function, \eqref{eq:overall objective2}, represents cost of travelling plus cost of energy usage for the deployed EVs. The increase in number of vehicles deployed can be thought of as the price paid to satisfy an uncertain demand.

\begin{table}[ht]
\begin{center}
\caption{Results Summary}
\begin{tabular}{l r r r r} 
 \hline
 \multicolumn{1}{c}{Code} & \multicolumn{1}{c}{Locations} & \multicolumn{1}{c}{Expected Passengers} & \multicolumn{1}{c}{Routes} &  \multicolumn{1}{c}{Objective \(\left[\$\right]\)}\\ 
 \hline\hline
    D-L5    & 5     & 224   & 2     & 1.32       \\     
    U-L5    & 5     & 224   & 2     & 1.32       \\     
    D-L10   & 10    & 492   & 4     & 0.87       \\ 
    U-L10   & 10    & 492   & 6     & 1.51       \\
    D-L15   & 15    & 696   & 5     & 1.12       \\
    U-L15   & 15    & 696   & 6     & 1.14       \\
    D-L20   & 20    & 904   & 8     & 1.33       \\
    U-L20   & 20    & 904   & 9     & 1.33       \\
    D-L25   & 25    & 1 196 & 11    & 1.10       \\
    U-L25   & 25    & 1 196 & 12    & 1.12       \\
    D-L30   & 30    & 1 396 & 12    & 1.31       \\
    U-L30   & 30    & 1 396 & 15    & 1.31       \\ 
    D-L50   & 50    & 2 144 & 18    & 1.43       \\
    U-L50   & 50    & 2 144 & 20    & 2.43       \\
    D-L75   & 75    & 3 208 & 35    & 2.12       \\
    U-L75   & 75    & 3 208 & 38    & 2.66       \\
    D-L100  & 100   & 4 424 & 65    & 2.50       \\
    U-L100  & 100   & 4 424 & 72    & 3.22       \\
    [1ex]
 \hline
\end{tabular}
\label{table:other results}
\end{center}
\vspace{-0.5cm}
\end{table}
\subsection{Energy Consumption}
\begin{figure}[ht]
     \centering
     \begin{subfigure}[b]{0.48\columnwidth}
         \centering
         \includegraphics[width=\columnwidth]{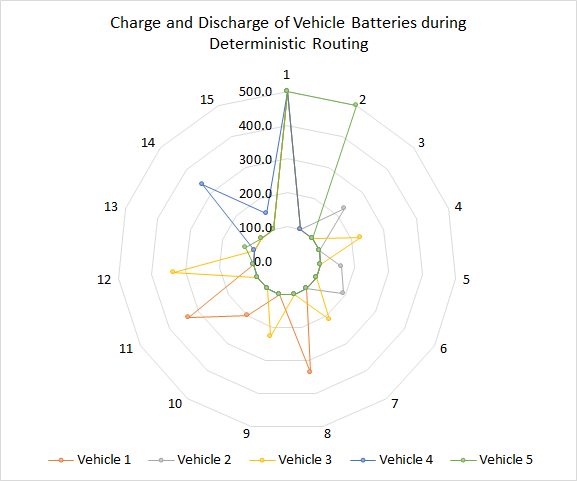}
         \caption{Energy levels resulting from  Deterministic model}
     \label{fig: Deterministic Routing}
     \end{subfigure}
     \hfill
     \begin{subfigure}[b]{0.48\columnwidth}
         \centering
         \includegraphics[width=\columnwidth]{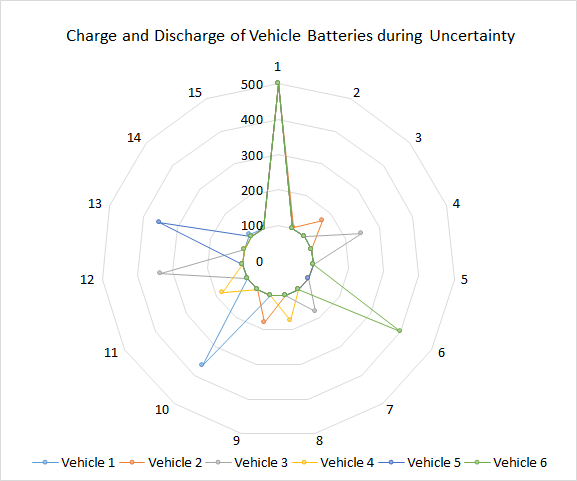}
         \caption{Energy levels  resulting from  chance-constrained model}
     \label{fig: Uncertain Routing}
     \end{subfigure}
        \caption{Energy levels during routing}
        \label{fig: Energy flow during Routing}
\end{figure}

Fig.~\ref{fig: Energy flow during Routing} presents the battery energy levels at different locations for scheduled routes for the 15-node case. The locations visited by each vehicle are represented in the spider web diagram while the height of spike is battery level of the EV.
Where the line representing each vehicle changes direction from downward-movement to upward-movement depicts a charging station.
In this 15-node problem, the charging stations are located at nodes 1, 2, 6, 8, 9 and 10. Take Vehicle 1 in Fig.~\ref{fig: Deterministic Routing} for example. Its route plan is 1-11-8-10-1 and this plan takes it through two charging stations at nodes 8 and 10. However, it only charges at one of them, node 8.
The route is complete when it returns to the depot and again recharges to the full energy level.
\subsection{Energy Consumption and $CO_{2}e$ Emissions}
In this section we analyse possible emissions' reductions when EVs are deployed instead of ICE vehicles for city-wide transportation. 
$5\%$ of power required by public transit in NYC is from electricity while $80\%$ comes from petroleum. For the types of vehicles listed, the emission factor is $2.26 KgCO_2e/unit$ of fuel and approximately 3.4 million gallons of gasoline and diesel fuel are used per day \cite{nycSIRR}. The average GHG emissions then stands at 19.23$kgCO_2e$ per vehicle.
On the other hand, the emission factor for electricity is $0.257 kgCO_2e/kWh$.

Tables~\ref{table:Emission for Deterministic Case} and \ref{table:Emission for Case with Chance Constraint} present the energy consumption and emissions results for the deterministic and the stochastic chance-constrained cases respectively. The value of GHG emissions for liquid fuel vehicles were estimated using the fuel emission factors above and a fuel-oil equivalent of energy used by EVs. Approximately, a 43\% reduction would be made in GHG emissions if EVs are utilized for shared taxi, as seen in Table~\ref{table:Emission for Deterministic Case}. When shared vans are used, which can accommodate between 8 and 10 passengers, we observe even greater savings. Van utilization lowers the number of vehicles, leading to a reduction in total $CO_2e$ emissions. If the energy consumption by vans and taxis is scaled based on number of passengers, emission by vans would be about 10\% of that of taxis.
\begin{table}[ht]
\begin{center}
\caption{Total Emissions for Deterministic Demand}
\begin{tabular}{l r r r r} 
 \hline
 \multicolumn{1}{c}{Type} & \multicolumn{1}{c}{Vehicles} & \multicolumn{1}{c}{Passengers}   & \multicolumn{1}{c}{Energy use} &  \multicolumn{1}{c}{Emissions} \\ 
 & & & \multicolumn{1}{c}{[kWh]} & \multicolumn{1}{c}{$\left[kgCO_2e\right]$} \\  
 \hline\hline
    Taxi (EV)           & 12    & 30    & 402.9 & 103.5      \\ 
    Van (EV)            & 5     & 150   & 228   & 58.6       \\ 
    Taxi (Liquid fuel)  & 12    & 30    & --     & 815.4     \\
    Van (Liquid fuel)   & 5     & 150   & --     & 461.04    \\
    [1ex]
 \hline
\end{tabular}
\label{table:Emission for Deterministic Case}
\end{center}
\vspace{-0.5cm}
\end{table}
\begin{table}[ht]
\begin{center}
\caption{Total Emissions with Stochastic Demand}
\begin{tabular}{l r r r r}
 \hline
 \multicolumn{1}{c}{Type} & \multicolumn{1}{c}{Vehicles} & \multicolumn{1}{c}{Expected}   & \multicolumn{1}{c}{Energy use} &  \multicolumn{1}{c}{Emissions} \\ 
 & & passengers & \multicolumn{1}{c}{[kWh]} & \multicolumn{1}{c}{$\left[kgCO_2e\right]$} \\ 
 \hline\hline
    Taxi (EV)           & 14 & 30  & 422.9 & 108.7      \\     
    Van (EV)            & 6  & 150 & 248.9 & 64.0       \\     
    Taxi (Liquid fuel)  & 14 & 30  & --     & 854.3      \\
    Van (Liquid fuel)   & 6  & 150 & --     & 503.75     \\
    [1ex]
 \hline
\end{tabular}
\label{table:Emission for Case with Chance Constraint}
\end{center}
\vspace{-0.3cm}
\end{table}

\section{Conclusion} \label{Conclusion}
This paper describes a new mixed integer linear programming formulation for EV routing which incorporates demand uncertainty, EV energy modelling and electric power grid utilization. Formulations were done such that bounds on confidence levels, simplicity and true representation of the systems we intended to describe were preserved. To reflect real life scenarios, we proposed a chance-constrained formulation for uncertainty. We thereafter based our case studies on inferences from real transport and electric network data of New York City. For our computational tests, differences between deterministic and stochastic model solutions were presented and analysed. This was done using the analysis of routing, energy consumption and emission of EVs deployed from the fleet.

Results from deterministic and stochastic models were contrasted for large number of locations and we noted that the increase in objective value for stochastic model is the extra cost of uncertainty. This increase in the objective value is reflected both in the amount of deployed vehicles and in their total energy consumption. Finally, we exhibit equivalent $CO_2$ emissions' reduction when EVs are used as taxis and as shared vehicles/vans. We show that there is a significant emissions' reduction when EVs are utilized for commercial transportation and we opine that this is a viable solution for climate change mitigation. The effect of a single aggregator of EVs on emissions reduction has been shown here. Future research direction could be investigating the cost/benefit for individual fleet owners if they collaborate with other fleet owners. Other advantages which include power system support and provision of ancillary services could be realized. 

\bibliographystyle{IEEEtran} 
\bibliography{refs.bib}

\end{document}